\documentclass{amsart}
\usepackage{amssymb}
\usepackage[dvips]{graphicx}
\usepackage{graphics}
\usepackage{psfrag}
\begin{document}
\newtheorem{thm1}{Theorem}[section]
\newtheorem{lem1}[thm1]{Lemma}
\newtheorem{rem1}[thm1]{Remark}
\newtheorem{def1}[thm1]{Definition}
\newtheorem{cor1}[thm1]{Corollary}
\newtheorem{defn1}[thm1]{Definition}
\newtheorem{prop1}[thm1]{Proposition}
\newtheorem{ex1}[thm1]{Example}
\newtheorem{alg1}[thm1]{Algorithm}


\title[toric ideals of graphs]{On complete intersection toric ideals of graphs}
\author{Christos Tatakis}
\address{Department of Mathematics, University of Ioannina,
Ioannina 45110, Greece }
\email{chtataki@cc.uoi.gr}
\author{Apostolos Thoma }
\address{Department of Mathematics, University of Ioannina,
Ioannina 45110, Greece }
\email{athoma@uoi.gr}
\thanks{}

%    \subjclass is required.
%    The 2010 edition of the Mathematics Subject Classification is
%    now available.  If you are citing a classification from the
%    new scheme, use the following input coding instead.
%\subjclass[2010]{Primary }

\date{}

\dedicatory{}

\begin{abstract}

\par We characterize the graphs $G$ for which their toric ideals $I_G$ are complete intersections.  
In particular we prove that for a connected graph $G$ such that $I_G$ is complete intersection
 all of its blocks are bipartite except of at most two.  We prove that toric ideals of graphs which are complete intersections are circuit ideals. 
The generators of the toric ideal correspond to even cycles of $G$ except of at most one generator, which corresponds to two edge
 disjoint odd cycles joint at a vertex or with a path.
 We prove that the blocks of the graph satisfy the odd cycle condition. Finally we characterize all complete intersection
toric ideals of graphs which are normal.
\end{abstract}
\maketitle

% \address[label2]{

\section{Introduction}

The complete intersection property of the toric ideals of graphs was first studied by L. Doering and T. Gunston in \cite{Doer}. In 1998 A. Simis proved 
that for a bipartite graph $G$ for which the toric ideal $I_G$ is complete intersection the number of chordless cycles of $G$ is equal
 to the number $m-n+r$, where  $m$ is the number of edges, $n$ the number of vetrices and $r$ the number of connected components of the graph $G$, see \cite{Simis}.
 Next year M. Katzman proved that for a bipartite graph $G$ the corresponding ideal $I_G$ is complete intersection if and only if any two chordless cycles 
have at most one edge in common, see \cite{Katzman}. Finally I. Gitler, E. Reyes,  and R. Villarreal determined  completely the form of the bipartite graphs 
for which the toric ideal $I_G$ is complete intersection. They are the {\em ring graphs}, see \cite{Gitler}. 
 Given a graph $H$, we call a path $P$ an $H$-path if $P$ is non-trivial 
and meets $H$ exactly in its ends. A graph $G$ is a {\em ring graph} if each block of $G$ which is not an edge or a vertex can
 be constructed from a cycle by successively adding $H$-paths of length at least two that meet graphs $H$ already constructed in two adjacent vertices.
 
\begin{thm1} \label{ring} [I. Gitler, E. Reyes,  and R. Villarreal \cite{Gitler}] If $G$ is a bipartite graph then $I_G$ is a complete intersection if and only if $G$
is a ring graph.
\end{thm1}

In this article we try to characterize complete intersection toric ideals of a general simple graph. 
Note that it is enough to answer the problem for a
connected graph, since for the toric ideal of a graph $G$ to be complete intersection it is enough that for every
connected component $G'$ of $G$ the ideal $I_{G'}$ to be complete intersection. 
In this article we will assume that all graphs considered are connected, except if stated otherwise.\\ The situation for a general graph is much more complicated than the case of a bipartite graph. For example, bipartite
complete intersection graphs are always planar, see \cite{Gitler},
but this is not the general case as the following example shows, see also \cite{Katzman}.

\begin{center}
\psfrag{A}{$e_{1}$}\psfrag{B}{$e_{2}$}\psfrag{C}{$e_{3}$}\psfrag{D}{$e_{4}$}\psfrag{E}{$e_{5}$}\psfrag{F}{$e_{6}$}\psfrag{G}{$e_{7}$}
\psfrag{H}{$e_{8}$}\psfrag{I}{$e_{9}$}\psfrag{J}{$e_{10}$}\psfrag{K}{$e_{11}$}
\includegraphics{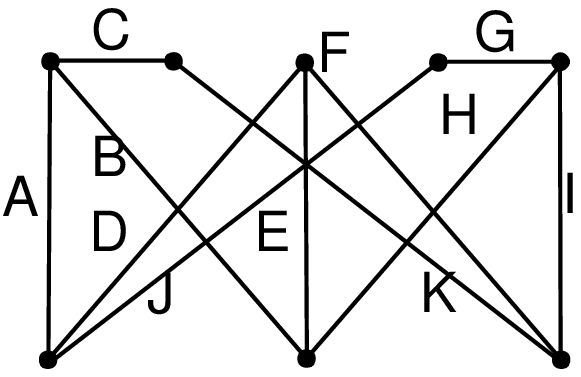}\\
{Figure 1}
%\caption{Figure 1}
\end{center}

Let $G$ be the graph with 11 edges and 8 vertices in Fig. 1. The height of the toric ideal $I_G$ is three, see \cite{Villa}, and  $I_G$ is generated by the
 binomials $e_1e_5-e_2e_4,e_5e_9-e_6e_8,e_3e_9e_{10}-e_1e_7e_{11}$ therefore
it is complete intersection. The graph $G$ is
 a subdivision of $K_{3,3}$ and therefore it is not planar, see \cite{planar}. A subdivision of a graph $G$ is any graph that can be obtained from $G$ by 
replacing edges by paths. 

Note also that the ideal of a general graph is much more complicated than the ideal of a bipartite graph. 
The generators of the toric ideal of a bipartite graph correspond to chordless even cycles of the graph. 
While the generators of the general graph have a more complicated structure, see Theorems \ref{primitive}, \ref{minimal}.
It is very interesting the fact that the generators of a complete intersection toric ideal are very simple, 
all of them correspond to even cycles with at most one exemption, see Theorem \ref{mostoneodd}.  Actually this is one of the properties that characterize 
complete intersection toric ideals of graphs, see Theorem \ref{last}.

In the second section 
 we review several notions from graph theory that will be usefull in the sequel. We define the toric ideal of a graph and we recall several results 
about the elements of the Graver basis, the circuits and the elements of a minimal system of generators of the toric ideal of the graph. 
The third section contain basic results about complete intersections toric ideals of graphs. The fourth section 
contains one of the main results of the article that
in a graph $G$ for which the toric ideal $I_G$ is complete intersection either all blocks are bipartite or all blocks are bipartite except one or 
all blocks are bipartite except two. In the case that there are exactly two non bipartite blocks they have a special position in the graph, the two blocks
 are {\em contiguous}. The fifth section contains the result that complete intersection toric ideals
 are circuit ideals and give a necessary and sufficient condition for a  graph $G$ to be complete intersection.
The final section proves that biconnected complete intersections graphs satisfy the odd cycle condition and  gives a necessary and sufficient condition for the edge
ring of a
complete intersection  graph to be normal.    

In the same problem, independently from us, I. Bermejo, I. Garc{\'i}a-Marco and E. Reyes  are working on providing
 combinatorial and alghorithmic  characterizations of general graphs such that their toric ideals are complete intersections in \cite{BerCI}.

\section{Toric Ideals of graphs}
  
  Let $A=\{{\bf a}_1,\ldots,{\bf a}_m\}\subseteq \mathbb{N}^n$
be a vector configuration in $\mathbb{Q}^n$ and
$\mathbb{N}A:=\{l_1{\bf a}_1+\cdots+l_m{\bf a}_m \ | \ l_i \in
\mathbb{N}\}$ the corresponding affine semigroup.  We grade the
polynomial ring $K[x_1,\ldots,x_m]$ over any field $K$ by the
semigroup $\mathbb{N}A$ setting $\deg_{A}(x_i)={\bf a}_i$ for
$i=1,\ldots,m$. For ${\bf u}=(u_1,\ldots,u_m) \in \mathbb{N}^m$,
we define the $A$-{\em degree} of the monomial ${\bf x}^{{\bf
u}}:=x_1^{u_1} \cdots x_m^{u_m}$ to be \[ \deg_{A}({\bf x}^{{\bf
u}}):=u_1{\bf a}_1+\cdots+u_m{\bf a}_m \in \mathbb{N}A.\]  The
{\em toric ideal} $I_{A}$ associated to $A$ is the prime ideal
generated by all the binomials ${\bf x}^{{\bf u}}- {\bf x}^{{\bf
v}}$ such that $\deg_{A}({\bf x}^{{\bf u}})=\deg_{A}({\bf x}^{{\bf
v}})$, see \cite{St}. For such binomials, we define $\deg_A({\bf
x}^{{\bf u}}- {\bf x}^{{\bf v}}):=\deg_{A}({\bf x}^{{\bf u}})$.

Let $G$ be a simple finite  connected graph on the vertex set $V(G)=\{v_{1},\ldots,v_{n}\}$ and let 
 $E(G)=\{e_{1},\ldots,e_{m}\}$ be the set of edges of $G$. We denote $K[e_{1},\ldots,e_{m}]$
 the polynomial ring in the $m$ variables $e_{1},\ldots,e_{m}$ over a field $K$.
 We
will associate each edge $e=\{v_{i},v_{j}\}\in E(G)$ with $a_{e}=v_{i}+v_{j}$ in the free
 abelian group generated by the vertices of $G$ and let $A_{G}=\{a_{e}\ | \ e\in E(G)\}$. We denote by $I_{G}$
 the toric ideal $I_{A_{G}}$ in
$K[e_{1},\ldots,e_{m}]$ and by $\deg_{G}$ the $\deg_{A_G}$. By $K[G]$ we denote the subalgebra of $K[v_1,\dots ,v_n]$ generated by all quadratic monomials
$v_iv_j$ such that $e=\{v_{i},v_{j}\}\in E(G)$. $K[G]$ is an affine semigroup ring and it is called the {\em edge ring} of $G$.  

A {\em cut vertex} (respectively cut edge) is a vertex (respectively edge) of the graph whose removal 
increases the number of connected components of the remaining subgraph. A graph is called {\em biconnected} if it is connected 
and does not contain  a cut vertex. A {\em block} is a maximal biconnected subgraph of a given graph $G$. 

A \emph{walk} of length $s$ connecting $v_{1}\in V(G)$ and $v_{s+1}\in V(G)$ is a finite
sequence of the form $$w=(\{v_{1},v_{2}\},\{v_{2},v_{3}\},\ldots,\{v_{q},v_{s+1}\})$$ with
each $e_j=\{v_{j},v_{j+1}\}\in E(G)$, $1\leq j\leq s$. An even (respectively odd) walk is a walk
of \emph{even} (respectively odd) length.
The walk $w$ 
is called \emph{closed} if $v_{s+1}=v_{1}$.
We call a walk $w'=(e_{j_1},\ldots,e_{j_t})$ a {\em subwalk} of $w$ if $e_{j_1}\cdots e_{j_t}|e_1\cdots e_s$.
A \emph{cycle} is a closed walk $$(\{v_{1},v_{2}\},\{v_{2},v_{3}\},\ldots,\{v_{s},v_{1}\})$$ with
$v_{i}\neq v_{j}$, for every $1\leq i < j \leq s$. For convenience by $\bf{w}$ we denote the subgraph of
$G$ with vertices the vertices of the walk and edges the edges of the walk $w$.
Given an even closed walk $$w=(e_{i_1},\ldots,e_{i_{2q-1}},e_{i_{2q}})$$ of the graph $G$ we denote by
$$E^+(w)=\prod_{k=1}^{q}e_{i_{2k-1}}={\bf e}^{w^+},\ E^-(w)=\prod_{k=1}^{q}e_{i_{2k}}={\bf e}^{w^-},$$ by $B_w$
the binomial $$B_w=\prod_{k=1}^{q}e_{i_{2k-1}}-\prod_{k=1}^{q}e_{i_{2k}}$$ belonging to the toric
ideal $I_G$, by $w^+, w^-$ the exponet vectors of the monomials $E^+(w), E^-(w)$ and by ${\bf w}^+, {\bf w}^-$ the sets $\{e_{i_1}, e_{i_3}, \dots , e_{i_{2q-1}}\}$,
$\{e_{i_2}, e_{i_4}, \dots , e_{i_{2q}}\}$ correspondigly.
Actually the toric ideal $I_G$ is generated by binomials of this form, see \cite{Villa}.
An even closed walk $w=(e_{i_1},\ldots,e_{i_{2q-1}},e_{i_{2q}})$ is said to be {\em primitive} if there exists no even closed subwalk
 $\xi$ of $w$ of smaller length such that $E^+(\xi)|E^+(w)$ and $E^-(\xi)|E^-(w)$.

Every even primitive walk $w=(e_{i_1},\ldots,e_{i_{2q}})$ partitions
the set of edges of $w$ in the two sets $\textbf{w}^+=\{e_{i_j}\ | \ j\ odd\}$ and $\textbf{w}^-=\{e_{i_j}\ | \ j \ even\}$,
otherwise if $e_{i_k}\in \textbf{w}^+\cap \textbf{w}^-$ then for the even closed subwalk $\xi=(e_{i_k},e_{i_k})$ we have $E^+(\xi)|E^+(w)$ and $E^-(\xi)|E^-(w)$. The edges of ${\bf w}^+$ are called odd edges of the walk
 and those of ${\bf w}^-$ are called even. {\em Sink} of a block $B$ of the graph ${\bf w}$ is a common vertex of two odd or two even edges of
the walk $w$ which belong to the block $B$.  

 Let $H$ be a subset of $V(G)$ and $G_H$ be the {\em induced graph} of $H$ in $G$, 
 which is the graph with vertices the elements of the set $H$ and edges the set of edges of $G$ where both vertices belong to $H$. For a given subgraph $F$  of $G$, 
an edge $f$
of the graph $G$ is called a \emph{chord} of the subgraph $F$  if the vertices of the edge
$f$ belong to $V(F)$ and $f\notin E(F)$. In other words an
edge is called chord of the the subgraph $F$ if it belongs to $E(G_{V(F)})$ but not in $E(F)$. A subgraph $F$ is called chordless if $F=G_{V(F)}$. For convenience by $G_w$ we denote the induced graph $G_{V({\bf w})}$, where $w$ is an even closed walk.

Let $w$ be an
even closed walk
$((v_{1},v_{2}),(v_{2},v_{3}),\ldots,(v_{2q},v_{1}))$ and
$f=\{v_{i},v_{j}\}$ a chord of $w$. Then $f$ \emph{breaks} $w$ in two
walks:
$$w_{1}=(e_{1},\ldots,e_{i-1}, f, e_{j},\ldots,e_{2q})$$ and
$$w_{2}=(e_{i},\ldots,e_{j-1},f),$$ where $e_{s}=\{v_{s},v_{s+1}\},\ 1\leq s\leq 2q$ and $e_{2q}=\{v_{2q},v_{1})\}.$
The two walks are both even or both odd. A chord $e=\{v_k,v_l\}$ is called {\em bridge} of a primitive walk $w$ if
there exist two different blocks $B_1,B_2$ of $\bf{w}$ such that $v_k\in B_1$ and $v_l\in B_2$. A chord
is called {\em even} (respectively {\em odd}) if it is not a bridge and breaks the walk in two even walks (respectively odd).
Thus we partition the set of chords of a primitive
even walk in three parts: bridges, even chords and odd chords.

\begin{def1}\label{creffec} Let
$w=(\{v_{i_{1}},v_{i_{2}}\}, \{v_{i_{2}},v_{i_{3}}\},\cdots ,
\{v_{i_{2q}},v_{i_{1}}\})$ be a primitive walk. Let
$f=\{v_{i_{s}},v_{i_{j}}\}$ and $f'=\{v_{i_{s'}},v_{i_{j'}}\}$ be two
odd chords (that means not bridges and $j-s,j'-s'$ are even) with $1\leq
s<j\leq 2q$ and $1\leq s'<j'\leq 2q$. We say that $f$ and $f'$
cross effectively in $w$ if $s'-s$ is odd (then necessarily $j-s',
j'-j, j'-s$ are odd) and  either $s<s'<j<j'$ or $s'<s<j'<j$.
\end{def1}

\begin{def1}\label{F_4}
We call an $F_4$ of the walk $w$ a cycle $(e, f, e', f')$ of
length four which consists of two edges $e,e'$ of the walk $w$
both odd or both even, and two odd chords  $f$ and $f'$ which
cross  effectively in $w$.
\end{def1}

A necessary and sufficient characterization of the primitive walks  of a graph, were given by
E. Reyes, Ch. Tatakis and A. Thoma in \cite[Theorem 3.2]{Thoma}:

\begin{thm1} \label{primitive} Let $G$ a  graph and $w$ an even closed walk of $G$. The walk $w$
is primitive if and only if
\begin{enumerate}
  \item every block of $\bf{w}$ is a cycle or a cut edge,
  \item every multiple edge of the walk $w$ is a double edge of the walk and a cut edge of $\bf{w}$,
  \item every cut vertex of $\bf{w}$ belongs to exactly two blocks and it is a sink of both.
\end{enumerate}
\end{thm1}

\begin{center}
 \includegraphics[scale=.6]{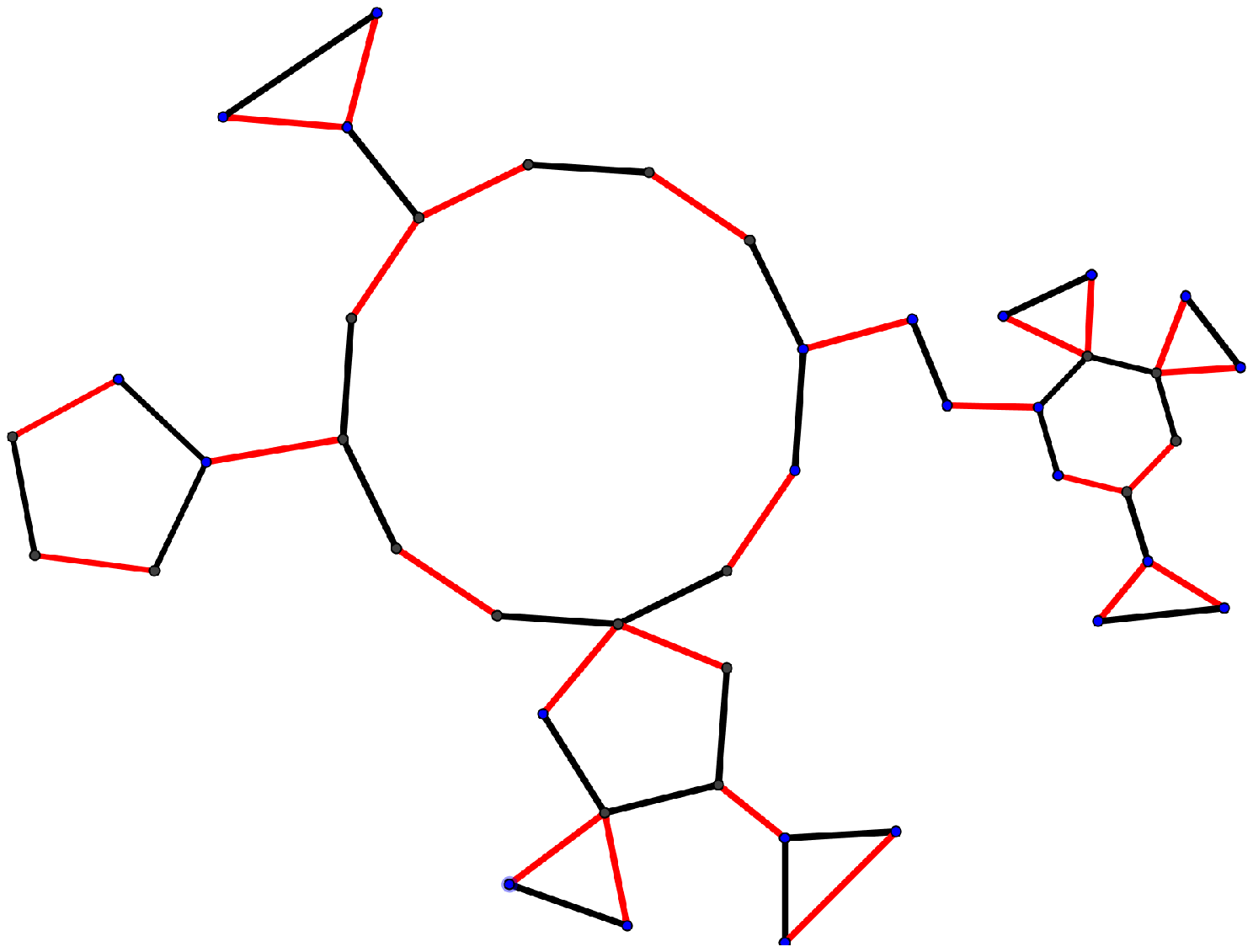}
\\Figure 2.
\end{center}

The following corollary were given by E. Reyes, Ch. Tatakis and A. Thoma in \cite[Corollary 3.3]{Thoma} and it describes the underlying graph of a primitive walk.

\begin{cor1} \label{primitive-graph}
Let $G$ be a  graph and let $W$ be a connected subgraph of $G$. The subgraph $W$ is the graph  ${\bf w}$ of a primitive walk $w$
 if and only if \begin{enumerate}
  \item  $W$ is an even cycle or
  \item  $W$ is not biconnected and
\begin{enumerate}
  \item every block of $W$ is a cycle or a cut edge and
  \item every cut vertex of $W$ belongs to exactly two blocks and separates the graph in two parts, the total number of edges
of the cyclic blocks in each part is odd.
\end{enumerate}
\end{enumerate}
\end{cor1}

In this case the walk $w$ passes through every edge of the cyclic blocks exactly once and from the cut edges twice.

A walk $w$ is primitive if and only if the binomial $B_w$ is primitive. The set of primitive binomials form the Graver basis of the toric ideal $I_G$. 
The Graver basis is important to us because every element of a minimal generating set of $I_G$ belongs to the Graver vasis of $I_G$, see \cite{St}.
We call {\em strongly primitive walk} a primitive walk that has not two sinks with distance one in any cyclic block, equivalently has not two adjacent cut vertices in any cyclic block.
For example the walk in Figure 1 is primitive but it is not strongly primitive, look for example at the cycle with six edges. We say that a binomial is minimal binomial
if it belongs to at least one minimal system of generators of $I_G$. 

The next theorem by E. Reyes, Ch. Tatakis and A. Thoma in \cite[Theorem 4.13]{Thoma} gives a necessary and sufficient characterization of minimal binomials of a 
toric ideal of a graph. This is the main theorem that made the results of this paper possible.

\begin{thm1} \label{minimal}
Let $w$ be an even closed walk. $B_{w}$ is a minimal binomial if
and only if \begin{enumerate}
  \item $w$ is strongly primitive,
  \item all the chords of $w$ are
odd and there are not two of them which cross  strongly effectively and
  \item no odd chord crosses an $F_4$ of the walk $w$.
\end{enumerate}
\end{thm1}

A necessary and sufficient characterization of circuits was given by R. Villarreal in \cite[Proposition 4.2]{Villa}:

\begin{thm1}\label{circuit}Let G be a  graph. The binomial $B\in I_G$ is circuit if and only if $B=B_w$, where $w$ is:
\begin{enumerate}
  \item an even cycle or
  \item two odd cycles intersecting in exactly one vertex or
  \item two vertex disjoint odd cycles joined by a path.
\end{enumerate}
\end{thm1}

\section{Complete intersection graphs}

The graph $G$ is called bipartite if it does not contain an odd cycle. The height of $I_{G}$ is equal to $h=m-n+1$ if $G$ is a bipartite graph or
$h=m-n$ if $G$ is a non-bipartite graph, where $m$ is the number
of edges of $G$ and $n$ is the number of its vertices, see \cite{Villa}. The toric ideal of $G$ is called a complete intersection if it can be generated by $h$ binomials.
We say that {\em a graph $G$ is complete intersection} if the ideal $I_G$ is complete intersection.

The problem of determining   complete intersection toric ideals has a long history starting with J.
Herzog in 1970 \cite{He} and finally solved by K. Fisher, W.~Morris and J. Shapiro in 1997 \cite{Fischer}.
For the history of this problem see the introduction of \cite{Morales}.

Next theorem says that the complete intersection property of a graph is hereditary property, in the sense that it holds also for
 all induced subgraphs.

\begin{thm1}\label{IGiff1}
The graph ${G}$ is complete intersection if and only if the graph $G_{H}$ is complete intersection
for every $H\subset V(G)$.
\end{thm1}
\noindent \textbf{Proof.} Let $B_{w_{1}},\ldots,B_{w_{s}}$ be a minimal system of generators of $I_{G_H}$, for some
even closed walks $w_i$ of $G$, $1\leq i\leq s$.
A minimal generator $B_w$ of $I_{G_H}$ is always
a minimal generator of $I_G$ since the property of being minimal generator
depends only on the induced graph $G_w$ of $w$ \ref{minimal}. Note that for a walk $w$ of $G_H$, the induced graph $G_w$ is the same in $G_H$ as in $G$. 
Therefore we can extend $B_{w_{1}},\ldots,B_{w_{s}}$ to a minimal system of generators $B_{w_1},B_{w_2},\ldots,B_{w_h}$ of $I_G$, $s\leq h$. 
The toric ideal $I_G$ is complete intersection therefore $B_{w_1},\ldots,B_{w_h}$ is a regular sequence. 
Since the ideal $I_{G}$ in $\mathbb{K}[e_1, \cdots, e_m]$ is homogeneous and non of the variables is a zero divisor in the edge ring $\mathbb{K}[G]=\mathbb{K}[e_1, \cdots, e_m]/I_G$, the sequence $B_{w_{1}},\ldots,B_{w_{s}}$
  is  regular  and therefore $I_{G_H}$ is a complete intersection toric ideal, see \cite{SSS}. \par
The converse is obvious since for  $H=V(G)$ we have $G=G_{H}$. \qed

The next proposition gives a very useful property of complete intersection toric ideals that will play a crucial role in the proofs of the theorems in the next sections.
\begin{prop1}\label{two}
  If  $G$ is complete intersection and $B_{w_{1}},\ldots,B_{w_{s}}$ is
a minimal set of generators of the ideal $I_{G}$ then there are no two walks $w_i$, $w_j$, $i\not=j$ such that
${\bf w}_i^+\cap {\bf w}_j^+\not= \emptyset$ and ${\bf w}_i^-\cap {\bf w}_j^-\not= \emptyset$,
or ${\bf w}_i^+\cap {\bf w}_j^-\not= \emptyset$ and ${\bf w}_i^-\cap {\bf w}_j^+\not= \emptyset$.
\end{prop1}
\noindent \textbf{Proof.} Let $B_{w_{1}}={\bf e}^{w_1^+}-{\bf e}^{w_1^-},\ldots,B_{w_{s}}={\bf e}^{w_s^+}-{\bf e}^{w_s^-}$ be
a minimal set of generators of the complete intersection toric ideal $I_G$. Then the matrix $M$ with rows $w_i^+-w_i^-$ is mixed dominating, see 
Corollary 2.10 \cite{FSh}. A matrix is called mixed if every row contains both a positive and a negative entry and dominating if it does not contain a square mixed submatrix.  Suppose that there exist  $B_{w_{1}},\ldots,B_{w_{s}}$
a minimal set of generators and two walks $w_i$, $w_j$, $i\not=j$ such that
${\bf w}_i^+\cap {\bf w}_j^+\not= \emptyset$ and ${\bf w}_i^-\cap {\bf w}_j^-\not= \emptyset$. Let $e_k\in {\bf w}_i^+\cap {\bf w}_j^+$ and 
$e_l\in {\bf w}_i^-\cap {\bf w}_j^-$. Then  the $2\times 2$ square submatrix taken from the $i, j$ rows and $k, l$ columns is mixed, 
contradicting the fact that $M$ is dominating. The proof of the other part is similar. \qed

It follows from the Proposition \ref{two} that if two edges are consequtive edges in two even closed walks $w_1$ and $w_2$
 in a complete intersection graph then both $B_{w_1}$, $B_{w_2}$ cannot belong in the same minimal system of generators of $I_G$. 
Also note that you cannot have in a minimal system of generators two circuits  with two odd cycles  and one of the cycles is the same in both, since any cycle contains at least three
edges and therefore there are at least two consequtive edges in common. For toric ideals of graphs Theroem \ref{minimal} determines the form a minimal binomial. 
Two minimal binomials sometimes belong to a minimal system of generators of the toric ideal, 
but for certain minimal binomials is impossible to find a minimal system of generators that contain both of them,
see \cite{ChKT}. For a toric ideal $I_A$ if two minimal binomials have different $A$-degrees then there exist a minimal system of generators for $I_A$ that contain both of them. 
But if they have the same $A$-degree sometimes there exist  a minimal system of generators for $I_A$ that contain both of them and some times not, for more details look at \cite{ChKT}.
For toric ideals of graphs the situation is simpler. Let $B_w$, $B_{w'}$ two minimal generators of $I_G$ then  there exist  
a minimal system of generators for $I_G$ that contain both of them if and only if $w$ and $w'$ are not $F_4$-{\em equivalent}. 

\begin{def1} Two primitive walks $w, w'$ differ by an
$F_4$, $\xi=(e_1,f_1,e_2,f_2)$, if $w=(w_1, e_1, w_2, e_2)$ and $w'=(w_1,
f_1,-w_2, f_2)$, where both $w_1, w_2$ are odd walks. Two primitive walks $w, w'$
are $F_4$-equivalent if either $w=w'$ or there exists a series of walks $w_1=w, w_2, \dots, w_{n-1}, w_n=w'$ such that
$w_i$ and $w_{i+1}$ differ by an $F_4$, where $1\leq i\leq n-1$.
\end{def1}
For more information about minimal system of generators of toric ideals of graphs see \cite{Thoma}.

\section{On the blocks of a complete intersection graph}

The Theorem \ref{oddblock} is one of the main results of the article and proves that if a complete intersection graph has $n$ blocks then at least $n-2$ of them are bipartite.
In the case that there are two nonbipartite blocks then they have to have a special position in the graph, they have to be {\em contiguous}.

\begin{def1} Two blocks of a graph $G$ are called contiguous if  there is a path
from the one to the other in which each edge of the path belongs to different block.
\end{def1}

Let $B(G)$ be the {\em block tree} of $G$, the bipartite graph with bipartition $(\mathbb{B},\mathbb{S})$ where $\mathbb{B}$ is the set of blocks of 
$G$ and $\mathbb{S}$ is the set of cut vertices of $G$,
$\{B, v\}$ is an edge if and only if $v\in B$. The leaves of the block tree are always blocks and are called {\em end blocks}. 
Let $B_k, B_i, B_l$ be blocks of a graph $G$. We call the block $B_i$ \emph{internal block}
 of $B_k,B_l$, if $B_i$ is an internal vertex in the unique path defined by $B_k,B_l$ in the tree $B(G)$.

\begin{thm1}\label{oddblock}
Let $G$ be a graph. If $G$ is complete intersection then either
\begin{enumerate}
\item all blocks of $G$ are
bipartite or 
\item all blocks are bipartite except one or
\item all blocks are bipartite except two which are contiguous.
\end{enumerate}
\end{thm1}
\noindent \textbf{Proof.}
Let $G$ be a  complete intersection graph
 and let $B_1,\ldots,B_t$ be the blocks of $G$. We assume that $G$ has  three
or more
non-bipartite blocks and let three of them be $B_m,B_k,B_l$. At least one of $B_m,B_k,B_l$ is not an internal
block of the other two, let it be $B_m$. We denote by $y_{i,j}$ the cut vertex of
$B_i$ which is the second vertex of the unique path $(B_i,\ldots,B_j)$ in the block tree $B(G)$, where $i,j\in \{m,k,l\}$.
Also we denote by $c_{i,j}$ an odd cycle of the block $B_i$ which contains the vertex
$y_{i,j}$ and with the smallest number of edges. Note that there exists at least one, since $B_i$ is a block which is non-bipartite, $i\in \{m,k,l\}$ .
 Let  $w_{m,k}=(c_{m,k},p_{m,k},c_{k,m},-p_{m,k})$, $w_{m,l}=(c_{m,l},p_{m,l},c_{l,m},-p_{m,l})$, where $p_{m,k},p_{m,l}$
are chordless paths from $y_{m,k}$ to $y_{k,m}$ and from $y_{m,l}$ to $y_{l,m}$ correspondingly and
 we can choose $c_{m,k}=c_{m,l}$ since $B_m$ is not an internal block of the other two. 
Note that whenever there is a path from a vertex to another then there is a chordless path between these two 
vertices. \\
We claim that the binomials $B_{w_{m,k}}$ and $B_{w_{m,l}}$ are minimal.
First $B_{w_{m,k}}$ is a circuit, see Theorem \ref{circuit} and therefore $w_{m,k}$ is primitive,
actually strongly primitive, see Theorem \ref{minimal}.
Note that $w_{m,k}$ has no bridges, since bridges are chords of the walk $w_{m,k}$
that their vertices are in different blocks of ${\bf w}_{m,k}$ which is impossible since:
a) $p_{m,k}$ is chordless, thus there is no bridge from the blocks of the path to themselves, b) $y_{m,k}$, $y_{k,m}$ are cut vertices, 
thus there is no bridge from the cycles to the path, and c) the odd cycles are of minimum length, therefore there is
no chord of the cycles incident to  $y_{m,k}$ or $y_{k,m}$. Also $w_{m,k}$ has no even chords since $c_{m,k}$, $c_{k,m}$ are odd cycles
 of minimum length.
 So all the chords of $w_{m,k}$ are odd.  Note that the odd chords of  $w_{m,k}$ are chords of either the cycle
$c_{m,k}$ or $c_{k,m}$.
There are not two of them which cross effectively, except if they form an $F_4$, otherwise
there will be an other odd cycle with strictly smaller number of edges than either $c_{m,k}$ or $c_{k,m}$ which passes from $y_{m,k}$ or $y_{k,m}$.
 Therefore by Theorem \ref{minimal} $B_{w_{m,k}}$ is minimal.
Similarly for the binomial  $B_{w_{m,l}}$. Note that $\deg_G(B_{w_{m,k}})\not=\deg_G(B_{w_{m,l}})$ 
thus they may belong to the same minimal system of generators of $I_G$.
A contradiction to Proposition \ref{two} since the cycle $c_{m,k}=c_{m,l}$ is contained in both walks. So the graph $G$ has at most two non-bipartite blocks.\\
Suppose that we are in the case that $G$ has exactly two non-bipartite blocks
and let them be $B_1$ and $B_2$.
 We will prove that they are contiguous.
Suppose not, then there exists at least one block $B_t$ such that every path from $B_1$ to $B_2$
 has at least two edges in $B_t$. Let $y_{t,1}$ and $y_{t,2}$ be the cut vertices of $B_t$ which are also vertices of the unique path $(B_1,\ldots,B_2)$
  in the block-tree $B(G)$. Since $y_{t,1}$ and $y_{t,2}$ belong in the same block $B_t$, 
there exist at least two internally disjoint paths of length at least two connecting them. Note that $\{y_{t,1}, y_{t,2}\}$ is not an edge of $G$, 
thus there 
are two different chordless paths from $y_{t,1}$ to $y_{t,2}$.
And so there exist at least two different chordless paths from $y_{1,2}$ and $y_{2,1}$.
Therefore by choosing the odd cycles $c_{1,2}$ and $c_{2,1}$ as the above construction, 
and  the two chordless paths $p_1, p_2$ from  $y_{1,2}$ and $y_{2,1}$ we get two even walks $w_1=(c_{1,2},p_1, c_{2,1}, -p_1)$,
$w_2=(c_{1,2},p_2, c_{2,1}, -p_2)$. 
As before,  there are no bridges in both $w_1, w_2$ and since all the chords of $c_{1,2}$ and $c_{2,1}$ are odd and there 
are not two of them which cross effectively (except if they form an $F_4$), each one of those paths will give a minimal generator of $I_G$.
Note that $\deg_G(B_{w_1})\not=\deg_G(B_{w_2})$ 
thus they may belong to the same minimal system of generators of $I_G$.
A contradiction to Proposition \ref{two} since the cycles $c_{1,2}$ and $c_{2,1}$ are contained in both walks.
 So $B_1$ and $B_2$ are contiguous blocks.
\qed

\section{Circuit ideals and complete intersections}

The first theorem of this section states an interesting property of toric ideals of graphs: complete intersection toric ideals of graphs are circuit ideals. 
Note that complete intersection toric ideals usually do not have this property. For more information on toric ideals generated by circuits see the  article \cite{Martinez} by
J. Martinez-Bernal and R. H. Villarreal and for circuit ideals  the article \cite{Bogart} by T. Bogart, A.N. Jensen and R.R. Thomas.

\begin{thm1}\label{circuit-ci}
Let $G$ be a  graph. If $G$ is a complete intersection
then every minimal generator of $I_{G}$ is a circuit.
\end{thm1}
\noindent \textbf{Proof.} Suppose that $I_G$ is a complete intersection toric ideal  that has
 a minimal generator $B_{w}$ which is not a circuit. Since the binomial $B_w$ is minimal
it is also primitive.  Therefore all of
the  blocks of ${\bf w}$ are cycles or cut edges, see Theorem \ref{minimal}.
Since $B_w$ is not a circuit it has at least three cyclic blocks from which at least two are odd, see Theorems \ref{minimal}
and \ref{circuit}.
The graph $G$ is complete intersection therefore the induced graph $G_w$ is complete intersection, from Theorem \ref{IGiff1},
where $G_w$ is the induced graph of ${\bf w}$ in $G$.
Note that the walk ${ w}$ has no bridges, since $B_w$ is a minimal generator, see Theorem \ref{minimal}, therefore there is a one to one correspondence
between the blocks of ${\bf w}$ and the blocks of $G_w$. Cut edges of ${\bf w}$ are cut edges of  $G_w$, but cyclic blocks of ${\bf w}$ may have chords
in $G_w$. At least two  of the blocks are non-bipartite, since they have an odd cycle.
Therefore from Theorem \ref{oddblock} there are exactly two. The end blocks of $B({\bf w})$ are always odd cycles of ${\bf w}$.  
Therefore the two non-bipartite blocks are the only end blocks of
 the block graph of $G_w$, which means that the block tree $B(G_w)$ is a path,
 see Corollary \ref{primitive-graph}. Let $B_1$, $B_2$ be the two odd cyclic blocks of ${\bf w}$ and $B_3$
be one of the other cyclic blocks, then $B_3$ will be an internal block of $B_1, B_2$. 
From Theorem \ref{oddblock} the two blocks $B_1$, $B_2$ are contiguous therefore
 there will be an edge of the block $B_3$ at the path which connects the two odd cycles
$B_1$, $B_2$ of $w$. If the edge belongs to the the walk $w$ then 
$w$ is not strongly primitive and if the edge does not belong to $w$ then it is a bridge of $w$, since its vertices are cut vertices of ${\bf w}$ and thus belong to two different blocks of ${\bf w}$. In both cases
 Theorem \ref{minimal} implies that $B_w$ is not a minimal generator, a contradiction.
Therefore ${\bf w}$ has at most two cyclic blocks and thus $B_w$ is a circuit, see Theorem \ref{primitive} and Theorem \ref{circuit}.
\qed

The next proposition will be usefull in the proof of Theorem \ref{noodd}.
\begin{prop1} \label{nochord} 
Let $G$ be a complete intersection graph and let $B_w$ be a minimal generator of $I_G$. If $w$ is not an even cycle then $w$ is chordless. 
\end{prop1}
\noindent \textbf{Proof.} From Theorem \ref{circuit-ci} the walk $w$  consists of two odd edge-disjoint cycles joint at vertex or with  a path, 
see Theorem \ref{circuit}. Thus $w$ is in the form $(c_1, p, c_2, -p)$, where $c_1, c_2$ are odd 
cycles $y_1, y_2$ are points of $c_1$ and $c_2$ correspondigly and $p$ is a path from $y_1$ to $y_2$ and it is possible that $y_1=y_2$ and $p$ to be empty. 
Since the binomial $B_w$ is minimal, the walk $w$ has no even chords and no bridges, see Theorem \ref{minimal}. Suppose that the walk $w$ had an odd chord $e=\{a,b\}$, then  from the definition of an odd chord both vertices belong to the same cycle.
Without loss of generality we can suppose
that both vertices $a, b$ belong to the cycle $c_1$. Then $c_1=(c_{11}, c_{12}, c_{13})$, where $ c_{11}, c_{12}, c_{13}$ are nonempty paths from $y_1$ to $a$, $a$ to $b$ and
$b$ to $y_1$, correspondigly. Among all possible such odd chords $e$ we choose the vertex $a$ in such a way that the length of $c_{11}$ is as small as possible. If there are more than one odd chord
with one vertex $a$ then we choose $b$ such that $c_{12}$ is as small as possible. By the choice of $a$ the walk $w_1=(c_{12}, \{b, a\}, -c_{11}, p, c_2, -p, c_{11})$ has no bridge from
$c_{12}$ to $c_{11}$. By the choice of $b$ there is no bridge from $c_{12} $ to the vertex $a$. Note that $w_1$ is not possible to have another bridge since $w$ has no bridges. Any chord of the odd cycle $(c_{12},e)$
is also a chord of $w$. $B_w$ is minimal generator of $I_G$ therefore if there exist such chords then all of them are odd chords of $w$ and not two of them  cross  strongly effectively and
   no odd chord crosses an $F_4$ of the walk $w$, see Theorem \ref{minimal}. Therefore also $B_{w_1}$ is minimal. 
Note that $\deg_G(B_{w})\not=\deg_G(B_{w_1})$ 
thus they may belong to the same minimal system of generators of $I_G$.
A contradiction to Proposition \ref{two} since the cycle $c_2$ is contained in both walks. \\
    Therefore
$w$ has no chord. \qed
\begin{thm1} \label{noodd}
Let $G$ be a biconnected complete intersection graph $G$. All  minimal generators of $I_{G}$ are in the form $B_w$ where $w$ is an even cycle.
\end{thm1}
\noindent \textbf{Proof.}
The Theorem \ref{circuit-ci} implies that all generators are circuits, thus to prove the theorem we will 
suppose that there is an even closed walk $w=(c_1,p,c_2,-p)$ of $G$ such that $B_w$ is a minimal generator of $I_{G}$ and we will arrive to a contradiction,
where $c_1,c_2$ are odd cycles of $G$,  $p=(v_1,\ldots,v_2)$ a path between them denoted by its vertices, 
$ V(c_1)\cap V(p)=\{v_1\}$ and $V(c_2)\cap V(p)=\{v_2\}$. 
Since $c_1,c_2$ are two edge-disjoint cycles of a biconnected graph, then there is at least 
one more path between them which is vertex disjoint from $p$. Let 
  $q=(x_1,y_1,\dots, y_2,x_2)$ be one of minimal length, where the vertex $x_1\in c_1$ and the vertex $x_2\in c_2$. 
Note that the length of $q$ is greater than one, since otherwise it will be just an edge which will be a bridge of $w$ 
and then $B_w$ will not  be a minimal generator of $I_G$. Note also that it may be $y_1=y_2$. Look at the graph induced by the graph 
${\bf w}\cup {\bf{q}}$. By Proposition \ref{nochord} $w$ has no chords and the  path $q$ has minimal length. 
Therefore the chords of ${\bf w}\cup {\bf{q}}$ are either edges from 
  the cycle $c_1$ to $y_1$, or from the cycle $c_2$ to $y_2$ and from $p$ to $q$
  except to the vertices $x_1, x_2$.  \\
  We claim that there are chords from the cycle $c_1$ to $y_1$. Suppose not. 
  Let $c_{11}$ be the path of greater length from $u_1$ to $x_1$ on the cycle $c_1$ and $c_{12}$
 be the path of smaller length from $u_1$ to $x_1$ on the circle $c_1$. Note that the cycle $c_1$ is odd and also the length of $c_{11}$ is
  greater than one. Denote by $c_{21}$ the path from $x_2$ to $u_2$ such that the cycle $w'=(c_{11}, \xi, c_{21}, -p)$ is even. 
There exist such path since the cycle $c_2$ is odd. 
  Consider $c$ to be the smallest even cycle in the form $(c_{11}, c')$ and the edges of it are edges or chords of $w'$. Note that there exist such 
cycle since $w'$   is in that form. Note that, $c$ is a cycle therefore it 
has no bridges and it is strongly primitive. Also from the minimality of the length of $c$ among even cycles of the form $(c_{11}, c')$, $c$ 
has no even chord and no two odd which cross strongly effectively 
and no odd that crosses an $F_4$ of $c$. Otherwise  the proofs of Propositions 4.8 and 4.12 of \cite{Thoma}
  show that 
  there exist two smaller even cycles and one of them is
  in the form $(c_{11}, c'')$, since there are no chords from $c_{11}$ (actually from $c$ to any vertex of $c'$). 
Then $B_c$ is a minimal generator of $I_G$.  Note that $\deg_G(B_{w})\not=\deg_G(B_{c})$ 
thus they may belong to the same minimal system of generators of $I_G$.
A contradiction to Proposition \ref{two} since the edges of $c_{11}$ are contained in both walks and length of $c_{11}$ is greater then one.\\
  Therefore there exist chords from the cycle $c_1$ to $y_1$ and similarly
  from the cycle $c_2$ to $y_2$.   Let $(u_{1,1}, u_{1,2},\dots, u_{1,s_1})$ be the cycle $c_1$ denoted by its vertices, where $u_{1,1}=u_1$. 
We consider all chords
   from the cycle $c_1$ to $y_1$ together with 
  the edge $\{x_1, y_1\}$ and denote them by 
  $e_1=\{y_1, u_{1,i_1}\}, \cdots, e_{t_1}=\{y_1, u_{1,i_{t_1}}\}$. Where $i_1<\dots<i_{t_1}$. Then the    
   cycles   $c_{1,j}=(y_1, u_{1,i_j},  u_{1,i_j+1},\dots,  u_{1,i_{j+1}})$, for $1\leq j\leq t_1-1$, 
and $c_{1,t_1}=(y_1, u_{1,i_{t_1}},  u_{1,i_{t_1}+1},\dots,  u_{1,s_1}, u_{1,1},\dots, u_{1,i_1})$ are chordless.
 If one of them was even then, since it does not have any chords, the corresponding binomial is a minimal generator
 of $I_G$, see Theorem \ref{minimal}. But it cannot have two consecutive edges which are in ${\bf w}^+$ or ${\bf w}^-$ therefore 
the only choice for the cycle is $(y_1,u_{1,s_1},u_{1,1}, u_{1,2})$. \\
  Similarly there cannot be two consecutive cycles $c_{1,j}, c_{1,j+1}$ odd since then 
the cycle $c=(y_1, u_{1,i_j},  u_{1,i_j+1},\dots,  u_{1,i_{j+1}}, u_{1,i_{j+1}+1}, u_{1,i_{j+2}} )$ 
is even with only one odd chord therefore the corresponding binomial is a minimal generator, see Theorem \ref{minimal}. 
 But $c$ cannot have two consecutive edges which are in ${\bf w}^+$ or ${\bf w}^-$ therefore the only choice for the resulting cycle 
is to be $(y_1,u_{1,s_1},u_{1,1}, u_{1,2})$, and the two original cycles where $(y_1,u_{1,1}, u_{1,2})$ and $(y_1,u_{1,s_1},u_{1,1})$.
  Since the cycle $c_1$ is even and each cycle $c_{1,j}$ consists of a part of $c_1$
 and two new edges, the number of odd cycles $c_{1,j}$ must be odd. And since among
 these cycles at most one can be even, then the number of odd cycles cannot be greater  than or equal to three since then you can 
find two consequtive cycles different from  $(y_1,u_{1,1}, u_{1,2})$ and $(y_1,u_{1,s_1},u_{1,1})$. 
 But then there is only one choice left, one cycle is even, the $(y_1,u_{1,s_1},u_{1,1}, u_{1,2})$, and one odd, 
the $(y_1, u_{1,2}, u_{1,3},\dots, u_{1,s_1},u_{1,1})$. 
   \\ A similar statement is true also for $y_2$  and the cycle $c_2$.
  The even closed walk $z=(y_1, u_{1,2}, u_{1,3},\dots, u_{1,s_1},u_{1,1}, y_1,\dots, y_2, u_{2,2}, u_{2,3},\xi', u_{2,s_2},u_{2,1}, y_2,-\xi',y_1)$ 
has no chords or bridges therefore $B_z$ is a minimal generator of $I_G$, see Theorem \ref{minimal}. 
But then from Proposition \ref{two} and that fact that $B_w$ is minimal generator we conlude that $s_1=3=s_2$.\\
Look at the graph ${\bf z}\cup {\bf w}$, the only chords of the graph can be from the path $p$ to the path $\xi '$.
  Let $e=\{a, b\}$ be the nearest chord to the vertex $u_{1,1}$, if there exist one,
  otherwise call $e$ the chord $\{ a=u_{2,1}, b=u_{2, 3}\}$.
  By the choice of the edge $e$ the cycle $o=(u_{1,1},\dots, a, b, \dots, y_1, u_{1,3}, u_{1,2})$ has no chord.\\
  There are two cases. First case: the cycle $o=(u_{1,1},\dots, a, b, \dots, y_1, u_{1,3}, u_{1,2})$ is even. 
Then $B_o$ is a minimal generator of $I_G$. But this is impossible since  the minimal generators $B_w$, $B_o$ have two consecutive edges 
in common,  $\{u_{1,3}, u_{1,2}\}$ and $\{u_{1,2}, u_{1,1}\}$, and $\deg_G(B_w)\not=\deg_G(B_o)$, see Proposition \ref{two}.
  \\ Second case: the cycle $o$ is odd, then the cycles $o_1=(u_{1,1},\dots, a, b, \dots, y_1,  u_{1,2})$, $o_2=(u_{1,1},\dots, a, b, \dots, y_1, u_{1,3})$ are both chordless and even. But then the $B_{o_1}, B_{o_2}$ are minimal generators of $I_G$.
  But this contradicts Proposition \ref{two} since the two cycles have all edges in common except two and $\deg_G(B_{o_1})\not=\deg_G(B_{o_2})$. \\
We conclude that all  minimal generators of $I_{G}$ are in the form $B_w$ where $w$ is an even cycle.\qed

\begin{thm1}\label{mostoneodd}
Let $G$ be a complete intersection graph. All minimal generators, except of at most one, of $I_G$ are in the form $B_w$ where $w$ is an even cycle. 
 The possible exceptional generator is a circuit whose two odd cycles belong to two different contiguous blocks.
\end{thm1}
\noindent \textbf{Proof.} In the case that all blocks of $G$ are bipartite or all exept one then there is no generator in the form $B_w$ 
where $w=(c_1,p,c_2,-p)$ 
with $c_1, c_2$ odd cycles, see Theorem \ref{noodd}. In the case that $G$ has two contiguous nonbipartite blocks then according to the proof of Theorem \ref{oddblock} there is one generator in the form $B_w$,
where $w=(c_1,p,c_2,-p)$ is an even closed walk, where $c_1,\ c_2$
are the unique odd chordless cycles of $B_1, B_2$ that are passing from  $y_{1,2}$ and $y_{2,1}$, correspondigly, and $p$ is the unique chordless path between them. 
Suppose that there is another generator in the form $B_{w'}$ where $w'$ is not an even cycle. Then from Theorem \ref{noodd} ${\bf w}'$ is not contained in the blocks $B_1$ or $B_2$. So ${\bf w}'$ consists of an odd chordless cycle
 $c_1'$ in the block $B_1$, see Proposition \ref{nochord}, an odd chordless cycle $c_2'$ in the block $B_2$ and a path $\xi$ from the one to the other. The path $\xi$
 consists of three paths. A chordless path $p_1$ from the cycle $c_1'$ to $y_{1,2}$, 
 the second is $p$ (since otherwise the path $\xi$ has a chord which plays the role
 of a bridge and destroys the minimality of the generator
 $B_{w'}$, see Theorem \ref{minimal}), and finally a  chordless path $p_2$ from the cycle $c_2'$ to $y_{2,1}$. Some of them may be empty. But then for the even closed 
 walk $w''=(c_1',p_1,p,c_2,-p,-p_1)$ we know that i) the two odd cycles $c_1', c_2$
 are chordless, ii) the path $(p_1, p)$ is chordless iii) there is no chord from 
 $c_1'$ to $p$, since $c_1'$ is in the block $B_1$, iv) there is no chord from 
 $c_2'$ to $p$, since $c_2'$ is in the block $B_2$ and v) there is no chord from
 $c_1'$ to $p_1$, since then it will be a bridge of $w'$ which is impossible from
 Theorem \ref{minimal}. Combining all these Theorem \ref{minimal} says that
  $B_{w''}$ is a minimal generator. The walks $w''$ and $w$ have more than two consecutive edges in common and $B_{w''}$, $B_w$ are minimal generators 
that they do not have the same $G$-degree, thus they may belong to the same minimal system of generators of $I_G$, contradicting Proposition \ref{two}. Thus there is no generator in the form $B_{w'}$ where $w'$ is not an even cycle.
\qed

For a block $B$ we denote by $I_B$ the toric ideal $I_G\cap K[e_i|e_i\in B]$, see \cite{St}.

The following result describes when a toric ideal $I_{G}$ is complete intersection.

\begin{thm1} \label{last}
Let $G$ be a  graph and let $B_{1},\ldots,B_{k}$ be its blocks.
$I_{G}$ is complete intersection toric ideal if and only if
\newline
i) all minimal generators, except of at most one, of $I_G$ are in the form $B_w$ where $w$ is an even cycle  and
\newline
ii) the ideals $I_{B_{i}}$ are complete intersection toric ideal for all $1\leq i\leq k$.
\end{thm1}
\noindent \textbf{Proof.}
Let $G$ be a  graph such that the toric ideal $I_{G}$ is complete intersection.
The first condition follows from Theorem \ref{mostoneodd} 
and the second from Theorem \ref{IGiff1} by choosing $H=V(B_i)$. \newline
Conversely, let $G$ be a graph and let $B_{1},\ldots,B_{k}$ be its blocks such that
 $I_{B_{i}}$ is complete intersection toric ideal for all $1\leq i\leq k$.
Note that every even cycle belongs to a unique block and all generators of the ideals $I_{B_{i}}$ 
correspond to even cycles, see Theorem \ref{noodd}. The number of minimal generators of the block $B_i$ is $m_i-n_i+1$ 
if $B_i$ is bipartite and $m_i-n_i$ if not. Therefore the total number of minimal generators of $I_G$ in the form $B_w$, where $w$ is an even cycle is
$$\sum _{i=1}^k (m_i-n_i+1) -j,$$ where $j$ is the number of nonbipartite blocks.
Note that $\sum _{i=1}^k m_i=m$, since every edge belongs to a unique block. $\sum _{i=1}^k n_i=n+\sum _{i=1}^c (deg(v_i)-1) $, where $v_i$ are cut vertices, 
 $deg(v_i)$ is the degree of $v_i$ as a vertex in the block tree $B(G)$ and $c$ is the number of cut vertices, since each cut vertex $v_i$ belongs
to $deg(v_i)$ blocks. $B(G)$ is a bipartite tree with bipartition $(\mathbb{B},\mathbb{S})$, where $\mathbb{B}$ is the set of blocks of $G$ and 
$\mathbb{S}$ is the set of cut vertices of $G$,
therefore   $\sum _{i=1}^c deg(v_i)$ is the number of edges of the tree $B(G)$ which is $k+c-1$. Combining all 
these we have that the total number of minimal generators of $I_G$ in the form $B_w$, where $w$ is an even cycle, is
$m-n-\sum _{i=1}^c deg(v_i)+c+k-j=m-n+1-j$.\\
We consider the following cases:\\
$j=0$, in this case the graph is bipartite and the total number of generators is $m-n+1$, since all minimal generators of $I_G$ in the form $B_w$, where $w$ is an even cycle.
Which means that the ${G}$ is a complete intersection.\\
$j=1$, in this case the graph is not bipartite and the total number of generators is $m-n$, since all minimal generators are in the form $B_w$, where $w$ is an even cycle.
Therefore ${G}$ is a complete intersection.\\
$j=2$, in this case the graph is not bipartite, thus its height is $m-n$ and the minimal generators in the form $B_w$, where $w$ is an even cycle, are $m-n-1$
so there must be exactly one more which is not in that form from condition $(i)$ and thus the total number of minimal generators is $m-n$ and  $G$ is
a complete intersection.\\
$j\geq 2$, in this case the graph is not bipartite, thus its height is $m-n$ and the minimal generators in the form $B_w$, where $w$ is an even cycle, 
are $m-n+1-j$
so there must be exactly one more which is not in that form from condition $(i)$ and thus the total number of minimal generators is $m-n+2-j$, which is less than the height, 
a contradiction to the generalized Krull's principal ideal theorem. \\ 
Therefore in all possible cases $G$ is a complete intersection.
\qed 

\section{The odd cycle condition and normality}

In this section we present  Theorems \ref{normal} and \ref{T_2} that are interesting  on their own, since they give us imformation about complete intersection graphs. 
But also they can be used to provide a necessary and sufficient condition for the edge ring of a complete intersection graph to be normal, see Theorem 
\ref{ci-normal}.
The normalization of the edge subring $K[G]$ was described explicitly by A. Simis,
W. V. Vasconcelos and R. V. Villarreal in \cite{SVV2} and by H. Ohsugi and T. Hibi in \cite{Hibi}. H. Ohsugi and T. Hibi related the normality of $K[G]$ with the {\em odd cycle condition}. 
\begin{def1} We say that a graph $G$ satisfies the odd cycle condition if for arbitrary two odd chordless cycles $c_1$ and $c_2$ in $G$, either $c_1, c_2$ have a common vertex or there exist an edge of $G$ joining a vertex of $c_1$ with a vertex of $c_2$. 
\end{def1}
For information about graphs satisfying the odd cycle condition see \cite{Fulkerson}, \cite{Ohsugi} and \cite{Stanley}. 
\begin{thm1} \label{odd cycle} [H. Ohsugi and T. Hibi \cite{Hibi}]  Let $G$ be a graph. Then the following conditions are equivalent:
\begin{itemize}
\item the edge ring $k[G]$ is normal,
\item the graph $G$ satisfies the odd cycle condition.
\end{itemize}
\end{thm1}

\begin{thm1} \label{normal}
Let $G$ be a biconnected complete intersection graph $G$. The graph $G$ satisfies the odd cycle condition and so the edge ring $K[G]$ is always normal.
\end{thm1}
\noindent \textbf{Proof.} Let $c_1$, $c_2$ be two chordless cycles of $G$ which have no
common vertex. We will prove that the subgraph ${\bf c}_1\cup {\bf c}_2$ has a chord. Suppose not. 
Let $p=(y_1, x_1, \dots, x_2, y_2)$ be the shortest path from $c_1$ to $c_2$, denoted by its vertices. 
The length of $p$ is greater than one, so it may be $x_1=x_2$.
The subgraph ${\bf c}_1\cup {\bf c}_2\cup {\bf p}$ has chords, otherwise the walk
$(c_1, p, c_2, -p)$ defines a minimal generator, which condraticts Theorem \ref{noodd}. 
Since $p$ is the shortest path it is chordless and there is no chord in ${\bf c}_1\cup {\bf c}_2$. Therefore all the chords of 
${\bf c}_1\cup {\bf c}_2\cup {\bf p}$ should be from the cycle $c_1$ to $x_1$ and from the cycle $c_2$ to $x_2$. 
 Let $(u_{1,1}, u_{1,2},\dots, u_{1,s_1})$ be the cycle $c_1$ denoted by its vertices, where $u_{1,1}=y_1$. We consider all chords
   from the cycle $c_1$ to $x_1$ together with 
  the edge $\{x_1, y_1\}$ and denote them by 
  $e_1=\{x_1, u_{1,i_1}\}, \cdots, e_{t_1}=\{x_1, u_{1,i_{t_1}}\}$. Where $1=i_1<\dots<i_{t_1}$. In the case that $t_1>1$,  the    
   cycles   $c_{1,j}=(x_1, u_{1,i_j},  u_{1,i_j+1},\dots,  u_{1,i_{j+1}})$, for $1\leq j\leq t_1-1$, 
and $c_{1,t_1}=(x_1, u_{1,i_{t_1}},  u_{1,i_{t_1}+1},\dots,  u_{1,s_1}, u_{1,1})$ are chordless and 
at least one of them is odd, say $c_{1,j}$, since $c_1$ is odd. Similarly if $t_2>1$ there must be an odd chordless cycle in the form $c_{2,k}$. \\
In the case that $t_1>1$ and $t_2>1$ let $w$ be the even closed path $(c_{1,j}, x_1, \dots, x_2, c_{2,k})$. \\
In the case that $t_1=1$ and $t_2>1$ let $w$ be the even closed path $(c_{1}, y_1, x_1, \dots, x_2, c_{2,k})$. \\
In the case that $t_1>1$ and $t_2=1$ let $w$ be the even closed path $(c_{1,j}, x_1, \dots, x_2, y_2, c_{2})$. \\
In the case that $t_1=1$ and $t_2=1$ let $w$ be the even closed path $(c_{1}, y_1, x_1, \dots, x_2, y_2, c_{2})$. \\
In all cases $w$ is chordless therefore by Theorem \ref{minimal} the binomial $B_w$ is minimal generator of $I_G$ contradicting Theorem \ref{noodd}.
We conclude that $G$ satisfies the odd cycle condition.
\qed
\begin{thm1}\label{T_2}
Let $G$ be a complete intersection graph such that it contains
two non-bipartite blocks $B_1,B_2$. 
Then each of the blocks $B_1,B_2$ contain atmost two odd chordless  cycles. Both  $B_1, B_2$ contain 
exactly an odd chordless cycle passing from the cut point $y_{1,2}$ and $y_{2,1}$ respectively. 
If any of them contained another one odd chordless  cycle then this cycle has distance one from the
 cut point $y_{1,2}$ if it is in $B_1$ or $y_{2,1}$ if it is in $B_2$. 
\end{thm1}

\textbf{Proof.}
Let $G$ be a complete intersection graph such that it contains
two non-bipartite blocks $B_1,B_2$. The two blocks are contiguous, see Theorem \ref{oddblock}. 
According to the proof of Theorem \ref{oddblock} there is one generator in the form $B_w$,
where $w=(c_1,p,c_2,-p)$ is an even closed walk, where $c_1, c_2$
are the unique odd chordless cycles of $B_1, B_2$ that are passing from  $y_{1,2}$ and $y_{2,1}$, 
correspondigly, and $p$ is the unique chordless path between them. 
Let $c_1=(\{y_{1,2},y_1\},\xi_1,\{y_2,y_{1,2}\})$, where $\xi_1$ is a path from $y_1$ to $y_2$. 
Suppose that the block $B_1$ contains another odd chordless cycle. Let $c$ be an odd chordless cycle different from $c_1$ 
and $\xi $ be a  path of smallest length from $c$ to $y_{1,2}$. 
Look at the even closed walk $w'=(c, \xi, p, c_2, -p, -\xi)$. Then $B_{w'}$ is not a minimal generator since it has a common cycle, the $c_2$, with $w$ and 
$\deg_G(B_{w'})\not=\deg_G(B_w)$. 
Since $c, c_2$ are chordless cycles, $\xi $ and $p$ are chordless paths and $c, \xi$
belong to the block $B_1$, $c_2$ belongs to the block $B_2$ 
and each edge of $p$ belongs to a different block there must be at least one chord (bridge of $w'$) from the cycle $c$ to the path $\xi$, 
see Theorem \ref{minimal}.  
And  since $\xi $ is a  path of smallest length from $c$ to $y_{1,2}$ any chord 
should be from $c$ to $x$ the second vertex of $\xi$. So, certainly $c$ does not passes from $y_{1,2}$ and this imply that $c_1$ is the only
 odd chordless cycle of the block $B_1$ that pass from $y_{1,2}$. 
Look at the induced graph of ${\bf c}\cup \{x\}$. Since $c$ is chordless any chord of ${\bf c}\cup \{x\}$ is from $c$ to $x$. 
 Let $(u_{1}, u_{2},\dots, u_{s})$ be the cycle $c$ denoted by its vertices. We consider all chords
   from the cycle $c$ to $x$ and denote them by 
  $e_1=\{x, u_{i_1}\}, \cdots, e_{t}=\{x_1, u_{i_t}\}$. Where $i_1<\dots<i_{t}$.  The    
   cycles   $o_{j}=(x, u_{i_j},  u_{i_j+1},\dots,  u_{i_{j+1}})$, for $1\leq j\leq t-1$, 
and $o_{t}=(x, u_{i_{t}},  u_{i_{t}+1},\dots,  u_{s}, u_{i_1})$ are chordless and 
at least one of them is odd since $c$ is odd. Without loss of generality we can suppose that it is $o_1$.  But then, the 
even closed walk $w''=(o_{1}, \xi', p, c_2, -p, -\xi')$ has no chords and bridges, where $\xi '$ is the subpath of $\xi $ from $x$ to $y_{1,2}$.
 Therefore by Theorem \ref{minimal}  $B_{w''}$ is a minimal generator  of $I_G$.
But this is not possible since  $w''$ has a common cycle, the $c_2$, with $w$, except if $B_{w''}=B_w$.   Therefore $o_{1}$
has to be $c_1$, $\xi'=\varnothing$, $u_{i_1}=y_1$, $u_{i_2}=y_2$ and $x=y_{1,2}$. 
Therefore $c$ has distance one from $y_{1,2}$ and it is in the form $(\xi_1, \xi_2)$, 
where $\xi_2$ is a path of even length from $y_2$ to $y_1$. The rest of the cycles $o_{i}$ are then necessary even, $2\leq i\leq t$. 
Note also that each of the cycles $o_i$ as well as $c$ are chordless. In case that $t>2$, the  cycle 
$(\{y_{1,2},u_{i_2}\},\{u_{i_2},u_{i_2-1} \}, \dots, \{u_{i_1+1},u_{i_1} \}, \xi_1,\{y_2,y_{1,2}\})$ is odd and chordless, pass from $y_{1,2}$ 
and is different from $c_1$. A contradiction, so $t=2$ and that means the subgraph ${\bf c}\cup \{x\}$ has only two chords, the  $\{y_{1,2},y_1\}$ and 
$\{y_2,y_{1,2}\}$.  \\It remains to prove that $c$ and $c_1$ are the only chordless odd
cycles in the block $B_1$. Suppose that there is another one $c'$. By repeating the proof we conclude that $c'$ is in the form $(\xi_1, \xi_3)$, 
where $\xi_3$ is a path of even length from $y_2$ to $y_1$ and the subgraph ${\bf c}'\cup \{x\}$ has only two chords, the  $\{y_{1,2},y_1\}$ and 
$\{y_2,y_{1,2}\}$. Then the cycles $w_1=(\{y_{1,2},y_1\},\xi_2,\{y_2,y_{1,2}\})$ and $w_2=(\{y_{1,2},y_1\},\xi_3,\{y_2,y_{1,2}\})$ are even and chordless
therefore $B_{w_1}$, $B_{w_2}$ belong to the same system of minimal generators of $I_G$ and share two consequtive edges. A contradiction. 
Therefore the block $B_1$ has at most two odd chordless cycles. 
\qed
  
\begin{def1} We say that a block is of type $T_i$ if it has $i$ chordless odd cycles.
 
\end{def1}
A bipartite block is of type $T_0$, while in a complete intersection graph 
with two non bipartite blocks the two blocks are either of type $T_1$ or of type $T_2$, 
from Theorem \ref{T_2}. Note that this is not true if the complete intersection graph has exactly one non bipartite block then it may be of
 higher type. For example the graph in Fig. 3 has type $T_4$ and it is complete intersection, since $I_G=(e_1e_3-e_2e_4, e_4e_6-e_5e_7, e_3e_5-e_8e_9)$ and 
and $h=9-6=3$. 

\begin{center}
 \includegraphics[scale=.4]{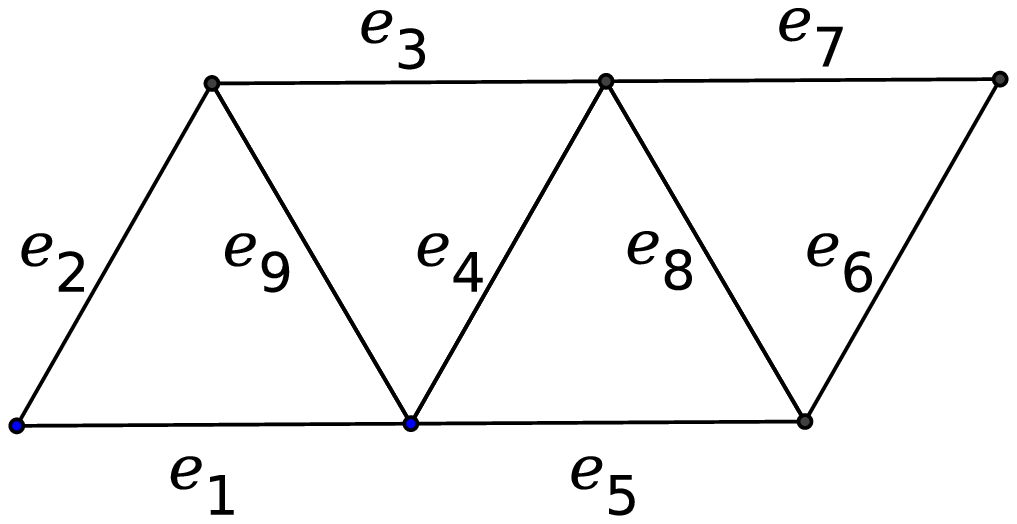}
\\Figure 3.
\end{center}

\begin{def1}
 
Two non-bipartite blocks are called strongly contiguous if
\begin{itemize}
\item  both are of type $T_1$ and they have distance at most one or
\item one is of type $T_1$ and the other of type $T_2$ and they have a common (cut) vertex.
\end{itemize}
\end{def1}
 It is easy to see that strongly contiguous blocks are always contiguous.

\begin{thm1} \label{ci-normal} Let $G$ be a complete intersection graph then $K[G]$ is normal if and only if $G$ has at most one non-bipartite
block or two which are  strongly contiguous.
\end{thm1}
\noindent \textbf{Proof.}  Let $G$ be a complete intersection graph such that $k[G]$ is normal, 
then $G$ satisfies the odd cycle property. From Theorem \ref{oddblock} we know that $G$ has at most one non-bipartite
block or two which are contiguous. In the first case we do not have anything to prove. Suppose that we are in the case that there 
are two contiguous blocks $B_1$, $B_2$. Then from Theorem \ref{T_2} they are either of type $T_1$ or of type $T_2$. \\In the case that
both are of type $T_1$ they have each exactly one odd chordless cycle passing from $y_{1,2}$ and $y_{2,1}$ respectively. 
Since $G$ satisfies the odd cycle property the two blocks $B_1, B_2$ have to  have distance at most one. \\In the case that 
one is of type $T_1$ and the other of type $T_2$, say $B_1$ is the first and $B_2$ the second, 
the block $B_1$ has only one odd chordless cycle passing through $y_{1,2}$, the block $B_2$ has two  odd chordless cycles, the one is passing from $y_{2,1}$ and the second has distance one from $y_{2,1}$. Since $G$ satisfies the odd cycle property the two blocks have to have  a common (cut) vertex, the $y_{1,2}=y_{2,1}$. Finally it is impossible to be both $B_1, B_2$ of type $T_2$,
since in this case there is an odd chordless cycle in $B_1$ with distance one from $y_{1,2}$ and there is an odd chordless cycle in $B_2$ with distance one from $y_{2,1}$. So these two odd cycles have distance at least two, contradicting the odd cycle property. So in all possible cases the two blocks are strongly contiguous.\\
For the converse, in the case that the graph $G$ has at most one non
 bipartite block then Theorem \ref{normal} implies that $G$ satisfies the odd cycle property and thus $K[G]$ is normal. In the case that $G$ has two  non-bipartite
blocks  which are  strongly contiguous Theorem \ref{T_2} implies that $G$ satisfies the odd cycle condition and thus $K[G]$ is normal. \qed

\end{document}